\documentclass{article}

\newcommand{\Om}{\Omega}

\newtheorem{theorem}{Theorem}[section]
\newtheorem{proposition}{Proposition}[section]
\newtheorem{lemma}{Lemma}[section]

\newtheorem{remark}{Remark}[section]
\newtheorem{remarks}{Remark}[section]
\newtheorem{definition}{Definition}[section]
\newcommand{\be}{\begin{equation}}
\newcommand{\ee}{\end{equation}}
\newcommand{\teta}{\theta}
\newcommand{\om}{\omega}

\newcommand{\ov}{\overline}
\newcommand{\wtilde}{\widetilde}
\newcommand{\R}{{\bf R}}

\newcommand{\sign}{{\rm sign}}
\renewcommand{\a }{\alpha }

\newcommand{\s }{\sigma }
\renewcommand{\d }{\delta }

\newcommand{\e }{\varepsilon }

\newcommand{\dps}{\displaystyle}

\begin{document}

\title{{\bf Bifurcation of free vibrations
for completely resonant wave equations
}}

\author{Massimiliano Berti, Philippe Bolle}
\date{}
\maketitle

{\bf Abstract:} 
We prove existence of small amplitude,   
$2\pi \slash \om$-periodic in time solutions of completely resonant 
nonlinear wave equations with Dirichlet boundary conditions
for any frequency $ \om $ 
belonging to a Cantor-like 
set of positive measure
and for a generic set of nonlinearities. 
The proof relies on a suitable Lyapunov-Schmidt decomposition 
and a variant of the Nash-Moser Implicit 
Function Theorem.
\\[2mm]
\footnote{Supported by 
M.I.U.R. Variational Methods and Nonlinear
Differential Equations.}
\\[1mm]
2000AMS subject classification: 35L05, 37K50, 58E05.

\section{Introduction and main result}

We outline in this note  recent results obtained in \cite{BB0} 
on the existence of small amplitude, $2\pi \slash \om$-periodic in time 
solutions
of the {\it completely resonant} nonlinear wave equation
\be\label{eq:main}
\cases{ u_{tt} - u_{xx} + f (x, u ) =  0 \cr
u ( t, 0 )= u( t,  \pi )  = 0}
\ee
where the nonlinearity $ f (x,u) = a_p (x) u^p + $
$O( u^{p+1} )$ with $p \geq 2$  
is analytic with respect to $u$ for $|u|$
small. More precisely, we assume 
\begin{itemize}
\item[{\bf (H)}] There is $ \rho > 0 $ such that $ \forall (x,u) \in (0,\pi)
  \times (-\rho,\rho)$, $f( x , u ) =  \sum_{k=p}^{\infty} a_k(x)
  u^k $, $ p \geq 2$, where $a_k \in H^1((0, \pi), \R)$ and 
$ \sum_{k=p}^{\infty} ||a_k||_{H^1} r^k < \infty$ for any $r\in
  (0,\rho)$. 
\end{itemize}
We look for periodic solutions of
(\ref{eq:main}) 
with frequency $ \om $ close to 
$1$ in a set 
of {\it positive measure}.
\\[1mm]
\indent 
Equation (\ref{eq:main}) is an infinite dimensional Hamiltonian
system
possessing an elliptic equilibrium at 
$ u = 0 $ with 
linear frequencies
of small oscillations $\om_j = j $, $ \forall j = 1, 2, \ldots $   
satisfying {\it infinitely many resonance} relations. 
Any solution $v = \sum_{j \geq 1} a_j \cos ( j t + \teta_j ) \sin (jx ) $ 
of the linearized equation at $u=0$, 
\be\label{eq:lin}
\cases{ u_{tt} - u_{xx} =  0 \cr
u ( t, 0 )= u( t,  \pi )  = 0}
\ee
is $2 \pi$-periodic in time.
For such reason equation (\ref{eq:main}) is called a {\it completely 
resonant} Hamiltonian PDE. 
\\[1mm]
\indent
Existence of periodic solutions of {\it finite} dimensional 
Hamiltonian systems close to a completely resonant elliptic equilibrium
has been proved by 
Weinstein, 
Moser 
and Fadell-Rabinowitz. 
The proofs are based on the classical 
Lyapunov-Schmidt
decomposition 
which splits the problem in two equations: the so called {\it range equation}, 
solved through the standard Implicit Function Theorem, and 
the {\it bifurcation equation} solved via variational arguments.  
\\[1mm]
\indent
For proving existence of small amplitude periodic solutions 
of 
completely resonant Hamiltonian PDEs like (\ref{eq:main}) 
two main difficulties must be overcome: 
\begin{itemize} 
\item[($i$)] a 
 ``{\it small denominators}'' problem which arises when solving the
 range equation; 
\item[($ii$)]  the  
presence of an {\it infinite dimensional} bifurcation equation:
which  solutions $ v $ of the linearized equation (\ref{eq:lin})
can be continued to solutions 
of the nonlinear equation (\ref{eq:main})? 
\end{itemize} 
The appearance of the small denominators problem $(i)$ is easily explained:
the  eigenvalues of the operator  $ \partial_{tt} -
\partial_{xx} $ 
in the space 
of functions  $u(t,x)$, $2 \pi / \om$-periodic in time 
and such that, say,  $u(t,.) \in  H^1_0 (0, \pi ) $ for all $t$, are 
$  - \om^2 l^2 + j^2 $, $l \in {\bf Z} $,
$ j \geq 1 $. 
Therefore, for almost every $ \om \in {\bf  R} $, 
the eigenvalues  accumulate to $ 0 $. As a consequence,  
for most  $\om$,  
the  inverse operator of $ \partial_{tt} -
\partial_{xx} $ is unbounded and the 
standard Implicit Function Theorem is not applicable.
\\[1mm]
\indent
The first existence results for small amplitude periodic solutions 
of (\ref{eq:main}) 
have been obtained 
in \cite{LS} 
for the specific nonlinearity $ f ( x, u ) = u^3 $
and 
periodic boundary conditions in $x$,
and in \cite{BP1} for $f(x,u)=u^3+ O(u^4)$, 
 imposing 
a ``strongly non-resonance'' condition on the frequency $\om$ 
satisfied in a {\it zero measure} set.
For such $ \om $'s the spectrum of $ \partial_{tt} - \partial_{xx} $ 
does not accumulate to 
$ 0 $ and so 
the small divisor problem ($i$) is bypassed.
The bifurcation equation (problem $(ii)$)  
is solved proving that, for $ f (x,u) = u^3 $, 
the $0^{th}$-order bifurcation equation possesses
{\it non-degenerate} periodic solutions.

In \cite{BB}-\cite{BB1}, for  
the same set 
of strongly non-resonant 
frequencies, existence and multiplicity of 
periodic solutions has been proved for {\it any } nonlinearity $ f(u) $. 
The novelty of \cite{BB}-\cite{BB1} was to solve the bifurcation equation 
via a variational principle
at fixed frequency which, jointly with 
min-max arguments, enables to find solutions of (\ref{eq:main})
as critical points of the Lagrangian action
functional. 

Unlike \cite{BP1}-\cite{BB}-\cite{BB1},
a new feature of the results we present in this Note
is that the set of frequencies $ \om $ for which we prove existence
of $2\pi \slash \om$-periodic in time solutions of (\ref{eq:main})
has positive measure.
\\[1mm]
\indent
Existence of periodic solutions for a positive measure set of frequencies 
has been proved in \cite{B2} in the case of periodic boundary conditions
in $x$ and for the specific nonlinearity 
$ f(x,u) = u^3 + \sum_{4 \leq j \leq d} a_j (x)u^j $
where the $a_j (x) $ are trigonometric cosine polynomials in $x$.
The nonlinear equation $u_{tt} - u_{xx} + u^3 =0$ with 
periodic boundary conditions 
possesses a continuum of small amplitude, 
analytic and non-degenerate 
periodic solutions 
in the form of travelling waves
$ u(t,x) = $ $ \d p_0 ( \om t + x) $.
With these properties at hand,
the small divisors problem ($i$) is solved in \cite{B2}
via a Nash-Moser Implicit function Theorem 
adapting the estimates of Craig-Wayne \cite{CW}. 

Recently, existence of periodic solutions of (\ref{eq:main}) for frequencies 
$\om$ in a positive measure set has been proved in 
\cite{GMP} using the Lindstedt series method
for odd analytic nonlinearities $ f ( u ) = a u^3 + O( u^5 )$
with $ a \neq 0 $.
The need for the dominant term $ a u^3 $ in the nonlinearity $f$ 
relies, as in \cite{BP1}, in the way the 
infinite dimensional bifurcation equation is solved. 
The reason for which $f(u)$ must be odd is that the solutions are
obtained as a sine-series in $ x $, see the comments
before Theorem \ref{thm:main}.
\\[1mm]
\indent
In \cite{BB0} we present a general method to prove existence 
of periodic solutions of the completely resonant wave equation
(\ref{eq:main}) with Dirichlet boundary conditions,
 not only for a positive measure set of 
frequencies $\om$, 
but also for a {\it generic} nonlinearity $ f(x,u) $ satisfying (H)
(we underline we do not require the oddness assumption $f(-x,-u) = -f(x,u) $),
see {\it Theorem \ref{thm:main}}.
\\[1mm]
\indent
Let's describe accurately our result.
Normalizing the period to $ 2 \pi $, we look for solutions 
$ u(t,x) $, $2 \pi$-periodic in time,  of the equation 
\be\label{eq:freq}
\cases{
\om^2 u_{tt} - u_{xx} + f (x, u ) = 0 \cr
u(t,0)= u(t, \pi) = 0
}
\ee
in  the real Hilbert space (which is actually 
a Banach algebra for $2s > 1$)
\begin{eqnarray*}
X_{\s,s} := \Big\{  u (t,x) = \sum_{l \in {\bf Z}} 
e^{i lt} \ u_l (x)   & \Big| &  u_l  \in H^1_0 ((0,\pi), {\bf C}), \ \ 
{\ov{u_l} (x)} = u_{-l} (x) \ \forall l \in {\bf Z}, \\ 
& & {\rm and} \  \ 
||u||_{\s,s}^2 := \sum_{l \in {\bf Z}} e^{2\s |l|} 
(l^{2s} + 1 ) ||u_l||^2_{H^1} < +\infty \Big\}. 
\end{eqnarray*}

For $ \s > 0 $ the space $ X_{\s,s}$                   
is the space of all $2\pi$-periodic in time functions with 
values in $ H^1_0 (( 0, \pi ), {\bf R}) $ which 
have a bounded analytic extension 
in the complex strip $|{\rm Im} \ t| < \s $ 
with trace function on $|{\rm Im } \ t| = \s $ belonging to  
$ H^s({\bf T}, H_0^1 ((0,\pi), {\bf C}))  $

The space of the solutions 
of the linear equation $ v_{tt} - v_{xx} = 0 $ that belong to 
$ X_{\s,s} $ is
$$
V := \Big\{  v (t,x) = \sum_{l \geq 1} 
\Big( e^{i lt} u_l + e^{- i lt} \ov{u_l} \Big)  
\sin (lx) \  \Big|  \ u_l \in {\bf C}   
\ {\rm and} \ || v ||_{\s,s}^2 = \sum_{l \in {\bf Z}} e^{2\s l} 
(l^{2s} + 1 )l^2 |u_l|^2 < +\infty \Big\}. 
$$
\indent

Let $\e:= \dps \frac{\om^2-1}{2}$.
Instead of looking for solutions of (\ref{eq:freq})
in a shrinking neighborhood of $0$ it is a
convenient devise to perform the rescaling $ u \to \d u $ 
with $ \d := |\e|^{1/p-1} $, obtaining
$$
\cases{
\om^2 u_{tt} - u_{xx} + \e g (\d, x, u ) = 0 \cr
u(t,0)= u(t, \pi) = 0}
$$
where 
$$
g(\d,x, u) :=  s^* \frac{f(x,\d u) }{\d^p} = s^* \Big( a_p(x) u^p + 
\d a_{p+1} (x) u^{p+1}
+ \ldots \ \Big)
$$
with $s^* :=\sign(\e )$, namely $ s^* = 1$ if $\om \geq 1$ and $s^*=-1$ if $\om < 1$. 
To fix the ideas, we shall consider  here  periodic solutions of 
frequency $\om >1$, so that $s^*=1$ and $\om=\sqrt{2 \delta^{p-1}+1}$.

If we try to
implement the usual Lyapunov-Schmidt reduction, 
i.e. to look for solutions
$ u = v + w $ with $ v \in V $ and $ w \in W: = V^\bot $,  
we are led to solve the bifurcation equation 
(sometimes called the ($Q$)-equation)
and the range equation (sometimes called the ($P$)-equation)
\be\label{eqs1}
\cases{
- \Delta v = \Pi_V g(\d, x,v + w) \  \ \ \qquad  \qquad (Q) \cr 
L_\om w = \e \Pi_W g(\d,x, v + w)   \ \qquad  \qquad (P)}
\ee
where 
$$ 
 \Delta v := v_{xx} + v_{tt}, \qquad 
\qquad   
L_\om := - \om^2 \partial_{tt} + \partial_{xx} 
$$ 
and 
$\Pi_V : X_{\s,s} \to V$,  $\Pi_W : X_{\s,s} \to W $
denote the projectors respectively on $V$ and $ W $.

Since $ V $ is infinite dimensional
a difficulty arises 
in the application of the method of \cite{CW} in presence
of small divisors :
 if $ v \in V \cap X_{\s_0,s} $ then the solution
$ w(\d , v) $ of the range equation, obtained with any Nash-Moser 
iteration scheme 
will have a lower regularity, e.g.  
$ w(\d, v) \in X_{\s_0 \slash 2,s} $. Therefore
in solving next the bifurcation equation for $v \in V $, 
the best estimate we can obtain is $v \in V \cap X_{\s_0 \slash 2,s+2} $,
which makes  the scheme incoherent. 
Moreover we have to ensure that
the $0^{th}$-order
bifurcation equation\footnote{We assume
for simplicity of exposition that the right hand side 
$ \Pi_V (a_p(x) v^p) $ is not identically equal to $0$ in $V$.
If not verified, the $0^{th}$-order non-trivial bifurcation equation 
will involve the higher 
order terms of the nonlinearity, see \cite{BB}. \label{note4}}, i.e.
the ($Q$)-equation for $ \d = 0 $, 
\be\label{eq:unpe}
- \Delta v = 
\Pi_V \Big( a_p(x) v^p \Big) 
\ee 
has solutions $ v \in V $ which are analytic, 
a necessary property to initiate an analytic Nash-Moser scheme
(in \cite{CW} this problem does not arise since,
dealing with {\it nonresonant} or {\it partially resonant}
Hamiltonian PDEs like
$ u_{tt} - u_{xx} + a_1 (x ) u = f(x,u) $,
the bifurcation equation is finite dimensional).
\\[1mm]
\indent
We overcome  this difficulty thanks to a reduction    
to a {\it finite dimensional} bifurcation equation (on a subspace  
of $V$ of dimension $N$ independent of $ \om $). This reduction can
be implemented,   
in spite of  the complete resonance of  equation (\ref{eq:main}),
thanks to the compactness of the operator $(-\Delta)^{-1}$.

We  introduce a  decomposition
$ V = V_1 \oplus V_2 $ 
where 
$$
\cases{
V_1 := \Big\{ v\in V  \ | \ v(t,x) = 
\sum_{l = 1}^N 
\Big( e^{i lt} u_l + e^{- i lt} \ov{u_l} \Big)
\ \sin (l x), \ u_l \in {\bf C} \Big\}  \cr 
V_2 := \Big\{ v \in V \ | \ v(t,x) = 
\sum_{l \geq N+1} \Big( e^{i lt} u_l + e^{- i lt}  \ov{u_l}\Big)
\sin (l x), \ u_l \in {\bf C} \Big\} 
}
$$
Setting $ v := v_1 + v_2 $,  
with 
$ v_1 \in V_1, v_2 \in V_2 $, (\ref{eqs1}) 
is equivalent to 
\be\label{eqs}
\cases{
- \Delta v_1 = \Pi_{V_1} g(\d, x,v_1 + v_2 + w) \ \qquad  \qquad (Q_1) \cr 
- \Delta v_2 = \Pi_{V_2} g(\d,x,v_1 + v_2 + w)  \ \qquad  \qquad (Q_2) \cr
L_\om w = \e \Pi_W g(\d,x, v_1 + v_2 + w)   \ \qquad  \qquad (P)}
\ee
where 
$\Pi_{V_i} : X_{\s,s} \to V_i$ ($ i = 1, 2 $), 
denote the orthogonal projectors 
on $V_i$ ($ i = 1, 2 $). 

Our strategy to find solutions of system (\ref{eqs}) is the
following. We solve first ({\it Step $1$})   
the $(Q_2)$-equation obtaining 
$ v_2 = v_2 (\d, v_1, w) \in V_2 \cap X_{\s , s }$
by a standard Implicit Function Theorem provided we have chosen  $N$ large enough 
and $\s $ small enough -depending on the nonlinearity $f$ 
but {\it independent of $\d$}.

Next ({\it Step $2$}) we solve the $(P)$-equation obtaining $w = w( \d, v_1) \in W 
\cap X_{\s \slash 2 ,s}$  
by means of a Nash-Moser Implicit Function Theorem
for $(\d , v_1) $ belonging to some Cantor-like set of parameters. 
A major role 
is played by the inversion of the {\it linearized operators}. 
Our approach -outlined in the next section- 
is much simpler than the ones usually employed 
and allows to deal 
nonlinearities which do {\sc not} satisfy the oddness assumption
$ f(-x,-u) = - f(x, u)$.
For this 
we develop $ u(t, \cdot ) \in H^1_0 (0,\pi)$ 
in time-Fourier expansion only.
Let us remark that $ H^1_0 (0, \pi) $ 
is the natural phase space
to deal with Dirichlet boundary conditions instead of the 
usually employed spaces $\{ u(x) = \sum_{j \geq 1 } u_j \sin (jx) \ | \
\sum_j e^{2 a j } j^{2 \rho } |u_j|^2 < + \infty  \}$, which 
force   the nonlinearity $ f $ to be odd.
We hope that the applicability of 
this technique  can go far beyond 
the present results.

Finally ({\it Step $3$})  
we solve the {\it finite dimensional} $(Q_1)$-equation for a generic set of 
nonlinearities obtaining $ v_1 = v_1 ( \d) \in V_1 $ 
for a set of $ \d $'s of positive measure.  

In conclusion we prove:

\begin{theorem}\label{thm:main} {\bf (\cite{BB0})}
Consider the completely resonant nonlinear wave 
equation (\ref{eq:main}) where the nonlinearity 
$f(x,u) =  a_p(x) u^p + O(u^{p+1}) $ satisfies assumption {\bf (H)}.


There exists an open and dense set ${\cal A}_p $ in $H^1((0,\pi), {\bf R})$ such that,
for all $a_p \in {\cal A}_p$, there is $ \s > 0 $ and a 
$C^\infty$-curve $[0, \d_0) \ni \delta \to u ( \delta ) \in X_{\sigma,s}$
with the following properties:
\begin{itemize}
\item  $(i)$ 
There exists $s^* \in \{-1,1\}$ and  a Cantor set ${\cal C}_{a_p}
\subset [0,\delta_0)$ satisfying 
\be\label{meas}
\lim_{\eta \to 0^{+}} \frac{{\rm meas} ({\cal C}_{a_p} 
\cap (0, \eta))}{\eta} = 1
\ee
such that, for all $ \delta \in
{\cal C}_{a_p}$, $ u ( \delta ) $ is a $2\pi / \om$-periodic in time 
solution of (\ref{eq:main}) with
$ \om = \sqrt{2 s^* \d^{p-1}+1}$;
\item  $(ii)$ 
$|| \wtilde{u}(\delta) - \delta u_0 ||_{\s,s } = O( \delta^2 )$ 
for some $ u_0 \in V \backslash \{0 \} \cap X_{\s , s }$
where $\wtilde{u}(\delta) (t,x)=u(\delta) (t/\om ,x)$. 
\end{itemize}

The conclusions of the 
theorem hold true for any nonlinearity 
$ f(x,u) = a_3 u^3 + \sum_{k\geq 4} a_k(x) u^k$, $a_3 \neq 0 $,
with $s^*= {\rm sign}(a_3)$. 
\end{theorem}


\section{Sketch of the proof}\label{sec:lya}

{\large {\bf Step 1: solution of the ($Q_2$)-equation.}}
The $0^{th}$-order bifurcation equation (\ref{eq:unpe}) 
is the Euler-Lagrange equation of the functional 
$\Phi_0 : V \to {\bf R}$ 
\be\label{func0}
\Phi_0 (v) = \frac{||v||_{H_1}^2}{2} - 
\int_\Om 
a_p(x) \frac{v^{p+1}}{p+1}
\ dx dt, \quad  \qquad \Om=(0,2\pi) \times (0,\pi). 
\ee
Assume for definiteness there is $ v \in V $ such that 
$ \int_\Om a_p(x) \frac{v^{p+1}}{p+1} > 0 $ 
(if  
the integral is $<0$ for some $v$, we can take $s^*=-1$ and 
 substitute $-a_p$ to $a_p$). 
%
Then $\Phi_0$  possesses 
by the Mountain-pass Theorem   
a non-trivial  critical set 
$ K_0 := \{ v \in V \ | \ \Phi_0'(v) =0, \Phi_0(v) = c \}$
which is compact for the $ H_1 $-topology, see \cite{BB}.
By a direct bootstrap argument 
any solution $ v \in K_0 $ of (\ref{eq:unpe}) belongs to 
$ H^k (V) $, $ \forall k \geq 0 $ and therefore is $ C^\infty $. 
In particular the Mountain-Pass solutions of (\ref{eq:unpe}) satisfy the 
{\it a-priori estimate} 
$
\sup_{v \in K_0} ||v||_{0,s+1} 
< R $ for some $ 0 < R < +\infty $.   
\\[1mm]
\indent
Solutions of the $(Q_2)$-equation are the 
fixed points of the nonlinear operator 
$ {\cal N}(\d,v_1,w, \cdot ) : V_2 \cap X_{\s,s} \to V_2 \cap X_{\s,s}$ 
defined by 
${\cal N}(\d, v_1, w, v_2 ) := $ 
$(- \Delta )^{-1} \Pi_{V_2}$ $g ( \d, x, v_1 + w + v_2 ).$
Using the {\it regularizing property} of
$ (- \Delta )^{-1} \Pi_2  $ 
we can prove that ${\cal N} $ is a 
contraction and then solve the ($Q_2$)-equation
in the space $ V_2 \cap X_{\s,s} $ for  $N$ large enough and for 
$0< \s < \ov{\sigma}$ 
($N$ and $\ov{\sigma}$ depend on $R$ but {\it not on $ \d $}).

\begin{lemma} \label{vdue}
{\bf (Solution of the ($Q_2$)-equation)}
There exist $ \ov{\s} > 0, N \in {\bf N}_+,\d_0>0$ such that,  
$\forall 0 < \s < \ov{\s} $, $\forall ||v_1||_{0,s+1} \leq 2R $,
$\forall  || w ||_{\s,s} \leq 1 $, $\forall |\d | \leq \d_0$, 
there exists a unique $v_2 = v_2(\d, w, v_1) \in X_{\s,s} $
with $||v_2 (\d, w,v_1)||_{\s,s} \leq 1 $
which solves the $(Q_2)$-equation.
Moreover $ v_2(\d, w, v_1) \in X_{\s,s+2} $.
\end{lemma}
 
Lemma \ref{vdue} implies, in particular, that any solution  
$ v \in K_0 $ of equation 
(\ref{eq:unpe}) is not only $ C^\infty $ but actually belongs 
to $X_{\s,s} $ and therefore is analytic
in $t$ (and hence in $ x $).  
\\[2mm]
{\large {\bf Step $2$: solution of the ($P$)-equation.}} 
By the previous step we are reduced to solve the $(P)$-equation 
with $ v_2 = v_2 (\d , v_1, w) $, namely
\be\label{eqP}
L_\om w = \e \Pi_W \Gamma (\d, v_1, w)
\ee
where $ \Gamma (\d, v_1, w)(t,x) := 
g(\d, x , v_1 (t,x) + w(t,x) + v_2(\d,v_1,w)(t,x) )$. 
\\[1mm]
\indent
The solution $ w = w(\d, v_1)$ of the $(P)$-equation 
(\ref{eqP}) is obtained 
by means of a Nash-Moser Implicit Function Theorem
for $( \d, v_1) $ belonging to a  Cantor-like set of parameters.
\\[1mm]
\indent
Consider the orthogonal splitting 
$ W = W^{(p)} \oplus W^{(p)\bot}$ where 
$W^{(p) }= $ $\{ w \in W \ | \ w = \sum_{l=0}^{ L_p} e^{ilt} \ w_l(x) 
\}$,  $W^{(p)\bot} =$ $ \{ w \in W \ | \ w = \sum_{l> L_p} 
e^{i lt} \ w_l(x) \}$
and $ L_p  = L_0 2^p $ for some large $L_0 \in {\bf N}$.
We denote by 
$P_p : W \to W^{(p)}$, 
$ P^{\bot}_p: W \to W^{(p)\bot} $
the orthogonal projectors
onto $W^{(p)}$, $W^{(p)\bot}$.
Define $ \s_0 := \ov{\s} $, 
the ``loss of analyticity at step $p$'' 
$ \gamma_p := \gamma_0 \slash (p^2+1)$ and  
$ \s_{p+1} = \s_p - \gamma_p $, 
$ \forall \ p \geq 0 $, 
with $ \gamma_0 > 0 $ small enough, such that the ``total
loss of analyticity'' $ \sum_{p\geq 0 } \gamma_p = $ $ \gamma_0 
\sum_{p \geq 0} 1 \slash (p^2+1) \leq$ $ \ov{\s}/ 2 $.

\begin{proposition} \label{Indu} 
{\bf (Nash-Moser iteration scheme)}
Let $ w_0 = 0 $ and
$ \  A_0 :=$ $\{ (\d,v_1) \ | \ |\d| < \d_0, $  $||v_1||_{0,s+1}$ $\leq 2R \}$.
There exist $\e_0 , L_0> 0 $ such that $\forall |\e | <\e_0 $,
 there exists a sequence 
$\{ w_p \}_{p \geq 0}$, $ w_p = w_p ( \d , v_1 ) \in W^{(p)} $, 
of solutions of
$$
L_\om w_p - \e P_p \Pi_W \Gamma (\d, v_1 , w_p ) = 0,
\leqno{(P_p)} 
$$
defined for 
$ ( \d, v_1 ) \in A_p \subseteq $ $ A_{p-1} \subseteq \ldots \subseteq 
A_1 \subseteq $ $ A_0 $.
For $ ( \d, v_1) \in A_\infty := \cap_{p \geq 0} A_p $, 
$ w_p (\d, v_1) $ totally 
converges in $ X_{\ov{\s}/2} $ to a solution $w (\d, v_1)$
of the $(P)$-equation (\ref{eqP})
with $|| w (\d, v_1) ||_{\ov{\s}/2,s} = O( \e ) .$

Moreover it is possible to define $ w ( \d, v_1 ) $
in a smooth way on the whole $ A_0 $:  
there exists a function $ \wtilde w (\d, v_1) \in C^\infty ( A_0, W) $
and a Cantor-like set $ B_\infty \subset A_\infty $
such that, if $ (\d, v_1) \in B_\infty \subset A_\infty $ 
then $ \wtilde{w}( \d, v_1 )$ solves
the ($P$)-equation (\ref{eqP}). 
\end{proposition}

Of course, the above proposition does not mean very much if we do not
specify $A_\infty$ or $B_\infty$. We refer to (\ref{spec}) for 
the definiton of $A_{p}$ and just say  that the set $B_\infty$ is sufficiently
large for our purpose. 
\\[1mm]
\indent
The real core of the Nash-Moser convergence proof 
-and where the analysis of the small divisors enters into play-
is the proof of the invertibility of the linearized operator
\begin{eqnarray*}
{\cal L}_p ( \d , v_1, w ) [ h ] & := &  
L_\om h - \e P_p \Pi_W D_w \Gamma ( \d , v_1 , w ) [ h ] \\
& = & 
L_\om h - \e P_p \Pi_W \Big( 
\partial_u g (\d ,x, v_1+w +v_2(\d, v_1,w) ) \Big[ h 
+ \partial_w v_2 (\d, v_1, w) [h] \Big] \Big), 
\end{eqnarray*}
where $w$ is the approximate solution obtained at a given stage
of the Nash-Moser iteration.
We do not follow the approach of \cite{CW} 
which is based on the Fr\"ohlich-Spencer techniques.

To invert ${\cal L}_p ( \d , v_1, w ) $, we distinguish a 
``diagonal  part'' $D$. 
Let 
$$
\cases{
a(t,x):=\partial_u g (\d ,x, v_1(t,x)+w(t,x)+v_2(v_1,w)(t,x)) \cr
a_0(x) := (1 \slash 2\pi )  \int_0^{2\pi} a(t,x) \ dt \cr
\ov{a}(t,x) := a(t,x)-a_0(x).}
$$
We can write 
$$
{\cal L}_p ( \d , v_1, w ) [ h ]
=  Dh - M_1h - M_2 h,
$$
where  $D $, $ M_1 $, $M_2 : W^{(p)} \to W^{(p)}$ are the linear operators  
\be\label{defin}
\cases{ 
Dh := L_\om h -\e P_p \Pi_W ( a_0 \  h) \cr
M_1 h := \e P_p \Pi_W (\ov{a}\  h)      \cr
M_2 h := \e P_p \Pi_W (a \  \partial_w v_2 [h]).}
\ee
We next diagonalize the operator $ D $ 
using Sturm-Liouville spectral theory.  
We find out that the  eigenvalues of $ D $ are $\om^2 k^2 - \lambda_{k,j} $,
$\forall  | k | \leq L_p $, $j\geq 1$, $ j \neq k $,    
and $\lambda_{k,j}$ satisfies the asympotic expansion 
\be\label{expa} 
\lambda_{k,j} = \lambda_{k,j}( \d, v_1, w) =
j^2 + \e M(\d, v_1,w)  + O \Big( 
\frac{\e ||a_0||_{H^1}}{j} \Big) \qquad
{\rm as } \qquad j \to + \infty,
\ee
where $M (\d, v_1,w) := (1\slash \pi )\int_0^\pi a_0(x) \ dx $. 

Assuming, for some 
$ \gamma > 0 $ and 
$ 1 < \tau < 2 $, the 
Diophantine condition (first order Melnikov condition)
\be \label{spec} \begin{array}{lll} 
( \d, v_1) \in A_p & := & \Big\{ (\delta, v_1) \in A_{p-1} \  \Big| \  
| \om k -j| \geq \dps \frac{\gamma}{(k+j)^{\tau}}, \ 
\Big| \om k -j -\e\dps  \frac{M(\delta , v_1 , w)}{2j} \Big| 
\geq \dps \frac{\gamma}{(k+j)^{\tau}}, \\ & &\\
& & \ \ \forall k \in {\bf N}, \ j \geq 1 \ {\rm s.t.}  \ k \neq j, 
\ \dps \frac{1}{3 |\e |} < k, \ j \leq L_p \Big\}
\subset A_{p-1}, 
\end{array} \ee  
all the eigenvalues of $D$ 
are {\it polynomially} bounded away from $0$, since 
$ \alpha_k := $
$ \min_{j \neq k, j \geq 1} |\om^2 k^2 - \lambda_{k,j}| \geq $ 
$ \gamma \slash k^{\tau - 1}$, $ \forall k $.
Therefore $ D $ is invertible
and $ D^{-1} $  
has sufficiently good estimates for the convergence of the
 Nash-Moser iteration. 

It remains to prove that the perturbative operators 
$ M_1 $, $ M_2 $ are small enough to get  the 
invertibility of the whole $ { \cal L}_p $.  
The smallness of $ M_2 $ 
is just a consequence of the regularizing 
property of $v_2: X_{\s,s} \to X_{\s,s+2} $ 
stated in Lemma \ref{vdue}.
The smallness of $ M_1 $ 
requires, on the contrary, 
an analysis of the ``{\it small divisors}'' $ \alpha_k $.
For our method it is sufficient  to prove that
$$
\a_k \a_l \geq c \gamma^2 |\e|^{\tau-1} > 0, \qquad \forall k \neq l \quad
{\rm with} \quad 
| k - l | \leq [\max\{ k, l \}]^{2-\tau \slash \tau}. 
$$
We underline again that this approach  works perfectly well for 
{\sc not} odd nonlinearities $ f $.
\\[2mm]
{\large {\bf Step $3$: solution of the ($Q_1$)-equation.}} 
Finally we have to solve the equation 
$$
- \Delta v_1 = \Pi_{V_1} {\cal G}(\d, v_1) \leqno{(Q_1)} 
$$
where $ {\cal G} ( \d, v_1 )(t,x) :=$ $ g ( 
\d, x, v_1(t,x) + \wtilde{w} ( \d, v_1 )(t,x)+
v_2 ( \d, v_1, \wtilde{w} ( \d, v_1 ))(t,x) ) $
and to ensure that there are solutions $( \d, v_1) \in B_\infty $  
for $ \delta $ in a set of positive measure (recall that 
if $(\delta,v_1) \in  B_\infty \subset A_\infty $,
then $ \wtilde{w} ( \d, v_1 ) $ 
solves the ($P$)-equation (\ref{eqP})).
Note that if 
$ \om = (1 + 2 \d^{p-1})^{1\slash 2} $ 
belongs to the zero measure set of ``strongly non-resonant'' frequencies 
used in \cite{BB}-\cite{BB1} then
$ ( \d, v_1 ) \in B_\infty, $ $\forall v_1 \in V_1 $
small enough. 
\\[1mm]
\indent
The finite dimensional $0^{th}$-order bifurcation equation, i.e. 
the ($Q_1$)-equation for $ \d = 0 $, 
$$
- \Delta v_1 = \Pi_{V_1} {\cal G}(0, v_1) = 
\Pi_{V_1} \Big( a_p(x) (v_1 + v_2 (0,v_1, 0))^p  \Big),
$$
is the Euler-Lagrange equation of the functional
$ {\wtilde \Phi}_0: V_1 \to {\bf R} $ where 
$ {\wtilde \Phi}_0 := \Phi_0 ( v_1 + v_2 (0,v_1, 0) )$
and $ \Phi_0 : V \to {\bf R}$ is the
functional  defined in (\ref{func0}).

It can be proved  
that if $a_p$ belongs to an  {\it open} and {\it  dense} subset ${\cal A}_p $ 
of $H^1 ((0, \pi), {\bf R})$, then 
$ \wtilde{\Phi}_0 : V_1 \to {\bf R}$ 
(or the functional that one
obtains when substituting $-a_p$ to $a_p$) 
possesses
a non-trivial  {\it non-degenerate} critical point
$ \ov{v}_1 \in V_1 $ 
and so, 
by the Implicit function Theorem, there exists 
a $C^\infty$-curve $ v_1 ( \cdot ) : \ [0, \d_0) \to V_1 $
of solutions of the $(Q_1)$-equation  
with $ v_1 (0) = \ov{v}_1 $.

The smoothness of  $ \d \to v_1 (\d)$ then implies that 
 $ \{ (\d, v_1 (\d)); \d >0 \} $ intersects 
 $ B_\infty $ in a set whose projection
on the $\d$ coordinate is the Cantor set
$ {\cal C}_{a_p} $ of Theorem \ref{thm:main}-$(i)$, 
 satisfying the measure estimate
(\ref{meas}). 
Finally $u(\d) = \d u_0 + O(\d^2)$ where 
$ u_0 := \ov{v}_1 + v_2 (0, \ov{v}_1 ,0) \in V $  
is a (non-degenerate, up to time translations) solution of 
the infinite dimensional bifurcation equation (\ref{eq:unpe}).

\noindent
Massimiliano Berti, SISSA, Via Beirut 2-4,
34014, Trieste, Italy, {\tt berti@sissa.it}. 
\\[1mm]
\noindent
Philippe Bolle,
D\'epartement de math\'ematiques, Universit\'e
d'Avignon, 33, rue Louis Pasteur, 84000 Avignon, France, 
{\tt philippe.bolle@univ-avignon.fr}.

\end{document}